\documentclass{article}

\usepackage{arxiv}

\usepackage[utf8]{inputenc} % allow utf-8 input
\usepackage[T1]{fontenc}    % use 8-bit T1 fonts
\usepackage{hyperref}       % hyperlinks
\usepackage{url}            % simple URL typesetting
\usepackage{booktabs}       % professional-quality tables
\usepackage{amsfonts}       % blackboard math symbols
\usepackage{amsmath}
\usepackage{nicefrac}       % compact symbols for 1/2, etc.
\usepackage{microtype}      % microtypography
\usepackage{graphicx}
\usepackage{natbib}
\usepackage{doi}

\title{The Relation Between the EVD or SVD of Summands and the EVD or SVD of the Sum}

\author{ \href{https://orcid.org/0000-0002-6345-0266}{\includegraphics[scale=0.06]{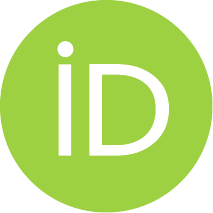}\hspace{1mm}Tsjerk A. Wassenaar} \\
\\
	Data Science for Life Sciences \\
	Hanze University of Applied Sciences \\
	Zernikeplein 11, 9747AS Groningen \\
    \\
    Molecular Dynamics Group \\
    University of Groningen \\
    Nijenborgh 7, 9747AG Groningen \\
    \\
    \texttt{t.a.wassenaar@rug.nl} \\
}

\renewcommand{\shorttitle}{The Relation Between the EVD or SVD of Summands and the EVD or SVD of the Sum}

\hypersetup{
pdftitle={The Relation Between the EVD or SVD of Summands and the EVD or SVD of the Sum},
pdfsubject={},
pdfauthor={Tsjerk A. Wassenaar},
pdfkeywords={Eigenvalue decomposition, EVD, Singular Value Decomposition, SVD, Principal Component Analysis, PCA},
}

\begin{document}
\maketitle

\begin{abstract}
In this work, we show how the eigenstructures of summands are related to that of the sum. In particular, we show that the sum of two positive semidefinite matrices can be written as the inner product of two block matrices $\mathbf{C} = \mathbf{A} + \mathbf{B} = \mathbf{PD}^2\mathbf{P}^T + \mathbf{QE}^2\mathbf{Q}^T = \begin{pmatrix} \mathbf{PD} & \mathbf{QE} \end{pmatrix}\begin{pmatrix} \mathbf{PD} & \mathbf{QE} \end{pmatrix}^T = \mathbf{Z}^T\mathbf{Z}$, such that the eigenvector decomposition of $\mathbf{C}$ can be obtained by projecting the eigenvectors from the block matrix product 

\[
\mathbf{Z}\mathbf{Z}^T = 
\begin{pmatrix} 
\mathbf{D}^2 & \mathbf{D}\mathbf{P}^T\mathbf{Q}\mathbf{E} \\
\mathbf{E}\mathbf{Q}^T\mathbf{P}\mathbf{D} & \mathbf{E}^2
\end{pmatrix}
\]

onto the block matrix $\mathbf{Z}^T = \begin{pmatrix} \mathbf{PD} & \mathbf{QE} \end{pmatrix}$. Next, it is shown that the result can be used to rewrite the SVD of a matrix inner product $\mathbf{X}^T\mathbf{Y}$, utilizing the eigenstructures of $\mathbf{X}$ and $\mathbf{Y}$. Finally, it is shown how this can be generalized to express the SVD of a sum of arbitrary matrices $\mathbf{H} = \mathbf{F} + \mathbf{G}$ in terms of the SVDs of the summands. The results may be useful in algorithms for eigendecomposition, specifically if the eigenproblem can be expressed in terms of a sum of matrices with simple eigenstructures as, e.g., in multi-omics research and generalized Procrustes analysis.
\end{abstract}

\keywords{Eigenvalue decomposition, EVD, Singular Value Decomposition, SVD, Principal Component Analysis, PCA}

\section{Introduction}

The eigenvalue decomposition (EVD) and the singular value decomposition (SVD) are cornerstone tools in mathematics, statistics, and physics. These decompositions 
underpin key methods in areas such as principal component analysis (PCA), signal processing, and machine learning. While the spectral properties of sums of matrices have 
been studied extensively, results have largely focused on bounds rather than structural relationships. Specifically, \citet{weyl1912} established inequalities for the 
eigenvalues of the sum of two positive semidefinite matrices in terms of the eigenvalues of the summands.

To date, no direct relationship between the eigenstructures (eigenvalues and eigenvectors) of the summands and those of the sum has been reported. In this work, we 
present such a relationship, offering insights into how the decompositions of the summands combine to form the decomposition of their sum. We also explore implications 
for the SVD of a matrix inner product ${\mathbf{X}}^T{\mathbf{Y}}$ and provide alternative formulations for its computation.

\section{Main Results}

\subsection{The eigenvector decomposition of the sum of two positive semidefinite matrices.}

\textbf{Theorem 1.} Let ${\mathbf{A}}$ and ${\mathbf{B}}$ be symmetric positive semidefinite square matrices of equal dimensions, with eigen decompositions ${\mathbf{PD}}^2{\mathbf{P}}^T$ and ${\mathbf{QE}}^2{\mathbf{Q}}^T$, respectively. Then the eigenvalues and eigenvectors of the sum ${\mathbf{C}} = {\mathbf{A}} + {\mathbf{B}}$ can be expressed in terms of the eigenstructure of a block matrix 

\begin{equation}
{\mathbf{Z}}^T{\mathbf{Z}} =
\begin{pmatrix} 
{\mathbf{D}}^2 & {\mathbf{D}}{\mathbf{P}}^T{\mathbf{Q}}{\mathbf{E}} \\
{\mathbf{E}}{\mathbf{Q}}^T{\mathbf{P}}{\mathbf{D}} & {\mathbf{E}}^2
\end{pmatrix}
\end{equation}

\textbf{Lemma 1.} For any matrix ${\mathbf{Z}}$, the products ${\mathbf{Z}}{\mathbf{Z}}^T$ and ${\mathbf{Z}}^T{\mathbf{Z}}$ share the same non-zero eigenvalues. 
Furthermore their eigenvectors are related as follows: the eigenvectors of ${\mathbf{Z}}^T{\mathbf{Z}}$ are projected onto ${\mathbf{Z}}$ to obtain those of 
${\mathbf{Z}}{\mathbf{Z}}^T$, and vice versa (scaled by ${\mathbf{Z}}^T$). This result can be shown, e.g., using the singular value decomposition (SVD).

\textbf{Proof} If ${\mathbf{A}} = {\mathbf{PD}}^2{\mathbf{P}}^T$ and ${\mathbf{B}} = {\mathbf{QE}}^2{\mathbf{Q}}^T$ then their sum can be written as the product 
${\mathbf{Z}}{\mathbf{Z}}^T$ of a block matrix ${\mathbf{Z}} = \begin{pmatrix} {\mathbf{P}}{\mathbf{D}} & {\mathbf{Q}}{\mathbf{E}} \end{pmatrix}$. Then using lemma 1 we 
can rewrite the objective of determining the eigenvectors of the sum of ${\mathbf{A}}$ and ${\mathbf{B}}$ as the objective of determining the eigenvectors of the matrix 
${\mathbf{Z}}^T{\mathbf{Z}}$ followed by projection of the resulting eigenvector on the matrix ${\mathbf{Z}}$.

\textbf{Corollary 1.} The number of non-zero eigenvalues and corresponding eigenvectors in ${\mathbf{P}}$ and ${\mathbf{Q}}$ is equal to the ranks of the matrices 
${\mathbf{A}}$ and ${\mathbf{B}}$, respectively. Thus, while the number of rows $n$ of ${\mathbf{Z}}$ is equal to that of ${\mathbf{A}}$, the number of columns $k$ is 
equal to the sum of the ranks of ${\mathbf{A}}$ and ${\mathbf{B}}$. Consequently, the block matrix ${\mathbf{Z}}^T{\mathbf{Z}}$ to diagonalize has dimensions $k$ by $k$ 
and may be considerably smaller than ${\mathbf{C}}$.

\textbf{Corollary 2.} The diagonal blocks ${\mathbf{D}}^2$ and ${\mathbf{E}}^2$ represent the eigenvalues of ${\mathbf{A}}$ and ${\mathbf{B}}$, respectively, and are thus diagonal. This can be leveraged to efficiently determine the eigenvectors of the block matrix.

\subsection{Expressing and solving the singular value decomposition as the sum of two matrices}

The result above can be used to rewrite and solve the SVD. A benefit of doing so, may be in the application of the SVD to the $(m \times n)$ product of two matrices ${\mathbf{X}}$ and ${\mathbf{Y}}$ of
dimensions $(k \times m)$ and $(k \times n)$, where $k < m + n$, which can be expressed as the EVD of a $(3k \times 3k)$ block matrix.

\textbf{Theorem 2.} Let ${\mathbf{X}}$ and ${\mathbf{Y}}$ be two matrices with dimensions $(k \times m)$ and $(k \times n)$, respectively. Then the SVD of the product 
${\mathbf{X}}^T{\mathbf{Y}}$ can be obtained from the EVD of the sum ${\mathbf{C}} = {\mathbf{A}} + {\mathbf{B}}$ of two positive semidefinite matrices with simple 
eigenstructures, where

\begin{equation}
{\mathbf{A}} = 
\begin{pmatrix} {\mathbf{X}}^T \\ {\mathbf{Y}}^T \end{pmatrix}
\begin{pmatrix} {\mathbf{X}} & {\mathbf{Y}} \end{pmatrix} =
{\mathbf{PD}}^2{\mathbf{P}}^T
\end{equation}

and

\begin{equation}
{\mathbf{B}} = 
\alpha {\mathbf{I}} - 
\begin{pmatrix}
{\mathbf{X}}^T{\mathbf{X}} & {\mathbf{0}} \\
{\mathbf{0}} & {\mathbf{Y}}^T{\mathbf{Y}}
\end{pmatrix} =
{\mathbf{QE}}^2{\mathbf{Q}}^T
\end{equation}

and solved according to \textbf{theorem 1} through the EVD on the block matrix

\begin{equation}
{\mathbf{Z}}^T{\mathbf{Z}} =
\begin{pmatrix} 
{\mathbf{D}}^2 & 
{\mathbf{D}}{\mathbf{P}}^T{\mathbf{Q}}_X{\mathbf{E}}_X &
{\mathbf{D}}{\mathbf{P}}^T{\mathbf{Q}}_Y{\mathbf{E}}_Y \\
{\mathbf{E}}_X{\mathbf{Q}}_X^T{\mathbf{P}}{\mathbf{D}} & 
{\mathbf{E}}_X^2 &
{\mathbf{0}} \\
{\mathbf{E}}_Y{\mathbf{Q}}_Y^T{\mathbf{P}}{\mathbf{D}} & 
{\mathbf{0}} &
{\mathbf{E}}_Y^2 
\end{pmatrix}
\end{equation}

where the subscripts $X$ and $Y$ refer to the corresponding blocks of ${\mathbf{B}}$

\textbf{Lemma 2.1.} The SVD of ${\mathbf{X}}^T{\mathbf{Y}}$ can be obtained from the EVD of the traceless augmented matrix.

\begin{equation}
\begin{pmatrix}
{\mathbf{0}} & {\mathbf{X}}^T{\mathbf{Y}} \\
{\mathbf{Y}}^T{\mathbf{X}} & {\mathbf{0}}
\end{pmatrix} =
\begin{pmatrix} {\mathbf{X}}^T \\ {\mathbf{Y}}^T \end{pmatrix}
\begin{pmatrix} {\mathbf{X}} & {\mathbf{Y}} \end{pmatrix}
-
\begin{pmatrix}
{\mathbf{X}}^T{\mathbf{X}} & {\mathbf{0}} \\
{\mathbf{0}} & {\mathbf{Y}}^T{\mathbf{Y}}
\end{pmatrix}
\end{equation}

\textbf{Lemma 2.2.} Adding a multiple of the identity matrix shifts the eigenvalues, while leaving the eigenvectors unaffected. This can be used to ensure that 
${\mathbf{B}}$ is positive semidefinite by choosing a scaling factor $\alpha$ equal to the smallest eigenvalue of the block diagonal matrix from lemma 2.1.

\textbf{Proof.} By lemma 2.1 and lemma 2.2 the SVD of the matrix product ${\mathbf{X}}^T{\mathbf{Y}}$ can be rewritten using the EVD of the sum of 2 positive semidefinite 
matrices ${\mathbf{A}}$ and ${\mathbf{B}}$, which can be solved in terms of the eigenstructures of these matrices by application of theorem 1.

\textbf{Note.} By lemma 2.2, all singular values are shifted by the factor $\alpha$.

\textbf{Corollary.} While the matrix product ${\mathbf{X}}^T{\mathbf{Y}}$ has dimensions $m \times n$, the matrix ${\mathbf{A}}$ has at most $\min(k, m + n)$ nonzero 
eigenvalues, while the matrix ${\mathbf{B}}$ has at most $\min(k, m) + \min(k, n)$. This means that the matrix ${\mathbf{Z}}^T{\mathbf{Z}}$ is a square matrix of 
dimension $\min(k, m + n) + \min(k, m) + \min(k, n)$. If $k \ll m + n$, this may result in a significant advantage for computing the eigenvalues and corresponding 
eigenvectors.

\subsection{Generalizing theorem 1 to the SVD of the sum of two arbitrary matrices}

Theorems 1 and 2 can be combined to provide a solution to rewrite the singular value decomposition of the sum of two matrices in terms of the decompositions of the summands.

\textbf{Theorem 3.} Let ${\mathbf{F}}$ and ${\mathbf{G}}$ be arbitrary matrices of equal size with singular value decompositions ${\mathbf{PD}}^2{\mathbf{Q}}^T$ and ${\mathbf{RE}}^2{\mathbf{S}}^T$, respectively. Then the singular values and singular vectors of the sum ${\mathbf{H}} = {\mathbf{F}} + {\mathbf{G}}$ can be expressed in terms of those of the summands.

\textbf{Proof.} The sum can be rewritten as the product ${\mathbf{H}} = {\mathbf{X}}^T{\mathbf{Y}}$ with 

\begin{equation}
{\mathbf{X}}^T = 
\begin{pmatrix} {\mathbf{PD}} & {\mathbf{RE}} 
\end{pmatrix}
\end{equation}

and

\begin{equation}
{\mathbf{Y}} = 
\begin{pmatrix} 
{\mathbf{D}}{\mathbf{Q}}^T \\ {\mathbf{E}}{\mathbf{S}}^T 
\end{pmatrix}
\end{equation}

By Theorem 2, the singular values and singular vectors of ${\mathbf{H}}$ can then be obtained through the eigenvalue decomposition of a derived block matrix, following from lemmas 2.1 and 2.2. Thus, the singular value decomposition of the sum can be expressed in terms of the singular value decompositions of the summands.

\textbf{Corollary 1.} Although this proof establishes the existence of a relation between the SVD of the sum and the SVD of the summands, the explicit nature of this relation is not immediately evident. To make it clearer, we can start from the augmented matrix from lemma 2.1 and rewrite this as a sum, using the square of the imaginary unit:

\begin{equation}
{\mathbf{C}} = \begin{pmatrix}
{\mathbf{0}} & {\mathbf{X}}^T{\mathbf{Y}} \\
{\mathbf{Y}}^T{\mathbf{X}} & {\mathbf{0}}
\end{pmatrix} =
\begin{pmatrix} {\mathbf{X}}^T \\ {\mathbf{Y}}^T \end{pmatrix}
\begin{pmatrix} {\mathbf{X}} & {\mathbf{Y}} \end{pmatrix}
+
i^2\begin{pmatrix}
{\mathbf{X}}^T & {\mathbf{0}} \\
{\mathbf{0}} & {\mathbf{Y}}^T
\end{pmatrix}
\begin{pmatrix}
{\mathbf{X}} & {\mathbf{0}} \\
{\mathbf{0}} & {\mathbf{Y}}
\end{pmatrix}
\end{equation}

This decomposition motivates rewriting the problem in terms of the matrix inner product ${\mathbf{Z}}{\mathbf{Z}}^T$ where

\begin{equation}
{\mathbf{Z}} = \begin{pmatrix}
{\mathbf{X}}^T &  i{\mathbf{X}}^T & {\mathbf{0}} \\
{\mathbf{Y}}^T & {\mathbf{0}} & i{\mathbf{Y}}^T
\end{pmatrix} = \begin{pmatrix}
{\mathbf{X}}^T & {\mathbf{0}} \\
{\mathbf{0}} & {\mathbf{Y}}^T
\end{pmatrix} \begin{pmatrix}
{\mathbf{I}} &  i{\mathbf{I}} & {\mathbf{0}} \\
{\mathbf{I}} & {\mathbf{0}} & i{\mathbf{I}}
\end{pmatrix} = {\mathbf{WJ}}
\end{equation}

\textbf{Corollary 2.} The eigenvectors of the augmented matrix ${\mathbf{C}} = {\mathbf{Z}}{\mathbf{Z}}^T = {\mathbf{WJ}}{\mathbf{J}}^T{\mathbf{W}}^T$ remain invariant under powers of ${\mathbf{C}}$. Thus, for some 
eigenvector $\begin{pmatrix} {\mathbf{u}} \\ {\mathbf{v}} \end{pmatrix}$ and eigenvalue $\lambda$, we can write:

\begin{equation}
\begin{split}
{\mathbf{C}}^p \begin{pmatrix} {\mathbf{u}} \\ {\mathbf{v}} \end{pmatrix} 
&= {\mathbf{Z}}({\mathbf{Z}}^T{\mathbf{Z}})^{p-1}{\mathbf{Z}}^T 
\begin{pmatrix} {\mathbf{u}} \\ {\mathbf{v}} \end{pmatrix} \\
&= {\mathbf{WJ}}({\mathbf{J}}^T{\mathbf{W}}^T{\mathbf{WJ}})^{p-1}{\mathbf{J}}^T{\mathbf{W}}^T
\begin{pmatrix} {\mathbf{u}} \\ {\mathbf{v}} \end{pmatrix} \\
&= {\mathbf{W}}({\mathbf{JJ}}^T{\mathbf{W}}^T{\mathbf{W}})^{p-1}{\mathbf{JJ}}^T{\mathbf{W}}^T
\begin{pmatrix} {\mathbf{u}} \\ {\mathbf{v}} \end{pmatrix} \\
&= \lambda^p \begin{pmatrix} {\mathbf{u}} \\ {\mathbf{v}} \end{pmatrix}
\end{split}
\label{eq:cor32}
\end{equation}

From this, it follows that the eigenvector, when projected onto the matrix ${\mathbf{JJ}}^T{\mathbf{W}}^T$, is an eigenvector of the matrix ${\mathbf{JJ}}^T{\mathbf{W}}^T{\mathbf{W}}$ and that projection of that vector onto ${\mathbf{W}}$ returns the eigenvector of ${\mathbf{C}}$. Although imaginary numbers are introduced in the decomposition, they are resolved in the rearrangement in equation \ref{eq:cor32}, resulting in only real-numbered intermediates.

\textbf{Corollary 3.} From corollary 2, it is evident that the eigenvectors can be obtained from the matrix ${\mathbf{JJ}}^T{\mathbf{W}}^T{\mathbf{W}}$. The first part of this product simplifies to

\begin{equation}
{\mathbf{JJ}}^T = \begin{pmatrix} {\mathbf{0}} & {\mathbf{I}} \\ {\mathbf{I}} & {\mathbf{0}} \end{pmatrix}
\end{equation}

The latter part of the product can be expanded in terms of the singular value decompositions of the original summands ${\mathbf{F}}$ and ${\mathbf{G}}$ to show the relations between those and the final result:

\begin{equation}
\begin{split}
  {\mathbf{W}}^T{\mathbf{W}}
  &= \begin{pmatrix}
    {\mathbf{X}}^T{\mathbf{X}} & {\mathbf{0}} \\
    {\mathbf{0}} & {\mathbf{Y}}^T{\mathbf{Y}}
  \end{pmatrix} \\
  &= \begin{pmatrix} \begin{array}{cc|cc}
      {\mathbf{D}}^2 & {\mathbf{DP}}^T{\mathbf{SE}} & {\mathbf{0}} & {\mathbf{0}} \\ 
      {\mathbf{DP}}^T{\mathbf{SE}} & {\mathbf{E}}^2 & {\mathbf{0}} & {\mathbf{0}} \\
      \hline
      {\mathbf{0}} & {\mathbf{0}} & {\mathbf{D}}^2 & {\mathbf{DQ}}^T{\mathbf{RE}} \\ 
      {\mathbf{0}} & {\mathbf{0}} & {\mathbf{ER}}^T{\mathbf{QD}} & {\mathbf{E}}^2 \\      
  \end{array} \end{pmatrix}
\end{split}
\label{eq:cor33}
\end{equation}

Taken together, this can be used to gain insight in the nature of the relation between the singular vectors of the summands and those of the sum: the eigenvector of the matrix ${\mathbf{W}}^T{\mathbf{W}}$ consists of two parts, one related to the product ${\mathbf{X}}^T{\mathbf{X}}$, capturing the relations between left singular vectors, and one related to ${\mathbf{Y}}^T{\mathbf{Y}}$, capturing relations between right singular vectors. Each of these parts itself thus consists of two subparts, one related to the first summand and one related to the other. The subparts are combined according to the alignment of the singular vectors of the summands, ${\mathbf{P}}^T{\mathbf{S}}$ and ${\mathbf{Q}}^T{\mathbf{R}}$, respectively. The multiplication by ${\mathbf{JJ}}^T$ then swaps the two parts of the eigenvector to capture the relation between the left and right singular vectors.

\textbf{Note.} The blocks in the matrix ${\mathbf{W}}^T{\mathbf{W}}$ corresponding to the left and the right singular vectors are each similar in structure to the block matrix in \textbf{theorem 1}.

\textbf{Corollary 4.} Corollaries 2 and 3 together suggest procedures for solving the singular value decomposition of a sum, given the decompositions of the summands using methods related to power iteration. This is particularly interesting, as the diagonal blocks are themselves diagonal, offering computational advantages.

\textbf{Corollary 5.} The major block diagonal structure of equation \ref{eq:cor33} together with the swapping of the two vector parts shows that both parts of the eigenvector can be obtained from one part, e.g., using the relation

\begin{equation}
\lambda^2 {\mathbf{w}} = \begin{pmatrix}
      {\mathbf{D}}^2 & {\mathbf{DP}}^T{\mathbf{SE}} \\ 
      {\mathbf{DP}}^T{\mathbf{SE}} & {\mathbf{E}}^2 
\end{pmatrix}
\begin{pmatrix}
  {\mathbf{D}}^2 & {\mathbf{DQ}}^T{\mathbf{RE}} \\ 
    {\mathbf{ER}}^T{\mathbf{QD}} & {\mathbf{E}}^2
\end{pmatrix}  {\mathbf{w}}
\end{equation}

\textbf{Corollary 6.} The approaches sketched above can also be applied to other decompositions, such as QR and LU, but this may not yield the sparsity, i.e., diagonal submatrices, obtained using eigen decompositions.

\textbf{Corollary 7.} The approaches sketched above can be trivially extended to sums of more than two parts, as these can likewise be written as inner products of block matrices.

\section{Conclusion}

In this work, we have demonstrated how the eigenvalue and singular value decompositions of a sum of matrices can be systematically related to the decompositions of the summands. These relationships are established by rewriting the sum as an inner product of block matrices, which reveals how the eigenstructures of the summands combine to form the eigenstructure of the sum. These insights suggest opportunities for more efficient computation, particularly for matrices with few rows and many columns, where smaller inner products can be utilized to derive the eigenvectors. Such matrices are prevalent in contemporary data science, especially in fields like omics studies, where numerous features are measured across relatively few instances. Beyond this, our results offer potential applications in other domains requiring large-scale matrix computations, providing a foundation for further exploration of efficient numerical methods and their extensions to other matrix factorizations.

\bibliographystyle{unsrtnat}
\bibliography{svdofsum}  

\end{document}